\documentclass{gtart_h}

\def\ifplaintex{\expandafter\ifx\csname documentclass\endcsname\relax}

\def\gtp{{\mathsurround=0pt\it $\cal G\mskip-2mu$eometry \&\ 
$\cal T\!\!$opology $\cal P\!$ublications}}  

\def\recd{{\small Received:\qua\receiveddate\ifx\reviseddate\relax
\else\qquad Revised:\qua\reviseddate\fi\par}} 

\def\lognumber#1{\def\thelognumber{#1}}
\def\volumenumber#1{\def\thevolumenumber{#1}}
\def\volumeyear#1{\def\thevolumeyear{#1}}
\def\papernumber#1{\def\thepapernumber{#1}}
\def\pagenumbers#1#2{\def\startpage{#1}\def\finishpage{#2}}
\def\published#1{\def\publishdate{#1}}

\def\received#1{\def\receiveddate{#1}}
\def\revised#1{\def\reviseddate{#1}}
\def\accepted#1{\def\accepteddate{#1}}
\def\asciititle#1{\def\theasciititle{#1}}

\def\asciiauthors#1{\def\theasciiauthors{#1}}
\def\asciiaddress#1{\def\theasciiaddress{#1}}
\def\asciiemail#1{\def\theasciiemail{#1}}

\long\def\asciiabstract#1{\long\def\theasciiabstract{#1}}
\def\asciikeywords#1{\def\theasciikeywords{#1}}

\let\\\par\let\thelognumber\relax\let\thevolumenumber\relax
\let\thepapernumber\relax\let\thevolumeyear\relax\let\startpage\relax
\let\finishpage\relax\let\publishdate\relax\let\receiveddate\relax
\let\reviseddate\relax\let\accepteddate\relax\let\theasciititle\relax
\let\theasciiauthors\relax\let\theasciiaddress\relax
\let\theasciiabstract\relax\let\theasciikeywords\relax

\let\theasciiemail\relax

\ifplaintex
\font\logobig=cmssbx10 scaled 3836
\font\logomed=cmssbx10 scaled 2557
\else
\font\logobig=cmssbx10 scaled 4200
\font\logomed=cmssbx10 scaled 2800
\fi

\long\def\makeagttitle{   
\count0=\startpage
\agt\hfill      
\hbox to 45truept{\vbox to 0pt{\vglue -13truept{\logomed A\kern -.37em{\logobig 
T}\kern -.38em G}\vss}\hss}
\break
{\small Volume \thevolumenumber\ (\thevolumeyear)
\startpage--\finishpage\nl
Published: \publishdate}

\vglue .25truein

{\parskip=0pt\leftskip 0pt plus
1fil\def\\{\par\smallskip}{\Large\bf\thetitle}\par\medskip} \vglue
0.05truein

{\parskip=0pt\leftskip 0pt plus 1fil\def\\{\par}{\sc\theauthors}
\par\medskip}%
 
\vglue 0.03truein 

{\small\leftskip 25truept\rightskip 25truept{\bf Abstract}\stdspace\theabstract

{\bf AMS Classification}\stdspace\theprimaryclass
\ifx\thesecondaryclass\relax\else; \thesecondaryclass\fi\par
{\bf Keywords}\stdspace \thekeywords\par}\vglue 7truept

}   

\ifplaintex
\font\phead=cmsl9 scaled 950
\font\pnum=cmbx10 scaled 913
\font\pfoot=cmsl9 scaled 950
\headline{\vbox to 0pt{\vskip -4.5mm\line{\small\phead\ifnum
\count0=\startpage ISSN 1472-2739 (on-line) 1472-2747 (printed)
\hfill {\pnum\folio}\else\ifodd\count0\def\\{ }%
\ifx\theshorttitle\relax\thetitle\else\theshorttitle\fi\hfill{\pnum\folio}
\else\def\\{ and }{\pnum\folio}\hfill\ifx\theshortauthors\relax\theauthors
\else\theshortauthors\fi\fi\fi}\vss}}
\footline{\vbox to 0pt{\vglue 0mm\line{\small\pfoot\ifnum\count0=\startpage
\copyright\ \gtp\hfill\else
\agt, Volume \thevolumenumber\ (\thevolumeyear)\hfill\fi}\vss}}
\else
\font=cmsl9 scaled 1050
\font\lnum=cmbx10 
\font=cmsl9 scaled 1050
\makeatletter
\def\@oddhead{{\small\lhead\ifnum\count0=\startpage ISSN 1472-2739 
(on-line) 1472-2747 (printed)\hfill {\lnum\number\count0}\else\ifodd\count0
\def\\{ }\ifx\theshorttitle\relax \thetitle \else\theshorttitle\fi\hfill
{\lnum\number\count0}\else\def\\{ and }{\lnum\number\count0}
\hfill\ifx\theshortauthors\relax 
\theauthors\else\theshortauthors\fi\fi\fi}}\def\@evenhead{\@oddhead}
\def\@oddfoot{\small\lfoot\ifnum\count0=\startpage\copyright\ \gtp\hfill\else
\agt, Volume \thevolumenumber\ (\thevolumeyear)\hfill\fi}
\def\@evenfoot{\@oddfoot}
\makeatother
\fi
\let\maketitlepage\makeagttitle
\let\makeshorttitle\maketitlepage
\let\maketitle\maketitlepage

\newwrite\gtoutfile
\long\gdef\makeheadfile{  
{\def\\{, }\def\s{ }
\immediate\openout\gtoutfile head.xxx
\immediate\write\gtoutfile{Proxy-for: \ifx\theasciiauthors\relax
\theauthors\else\theasciiauthors\fi\s<\ifx\theasciiemail\relax\theemail\else\theasciiemail\fi>}
\immediate\write\gtoutfile{\noexpand\\}
\immediate\write\gtoutfile{Authors: \ifx\theasciiauthors\relax
\theauthors\else\theasciiauthors\fi}
{\def\\{ }\immediate\write\gtoutfile{Title: \ifx\theasciititle\relax
\thetitle\else\theasciititle\fi}}
\immediate\write\gtoutfile{Subj-class: GT or SG, GR etc}
\immediate\write\gtoutfile{MSC-class: \theprimaryclass\ifx\thesecondaryclass\relax\else, \thesecondaryclass\fi}
\immediate\write\gtoutfile{Journal-ref: Algebr. Geom. Topol. \thevolumenumber\s
(\thevolumeyear) \startpage-\finishpage}
\immediate\write\gtoutfile{Comments: Published by Algebraic and
Geometric Topology at}
\immediate\write\gtoutfile{\s\s\s  http://www.maths.warwick.ac.uk/agt/AGTVol\thevolumenumber/agt-\thevolumenumber-\thepapernumber.abs.html}
\immediate\write\gtoutfile{\noexpand\\}
\immediate\write\gtoutfile{}
\ifx\theasciiabstract\relax
\immediate\write\gtoutfile{\theabstract}\else
\immediate\write\gtoutfile{\theasciiabstract}\fi
\immediate\write\gtoutfile{}
\immediate\write\gtoutfile{\noexpand\\}
\immediate\write\gtoutfile{}
\immediate\closeout\gtoutfile}}  

\def\maketitlepage{\makeagttitle\makeheadfile}
\let\makeshorttitle\maketitlepage
\let\maketitle\maketitlepage

\lognumber{18}
\volumenumber{5}
\volumeyear{2005}
\papernumber{18}
\pagenumbers{405}{418}
\received{21 January 2005} 
\revised{14 April 2005}
\accepted{28 April 2005}
\published{22 May 2005}

\usepackage{amssymb,amsmath}
\allowdisplaybreaks

\newtheorem{theorem}{Theorem}[section]
\newtheorem{lemma}[theorem]{Lemma}
\newtheorem{corollary}[theorem]{Corollary}
\newtheorem{conjecture}[theorem]{Conjecture}

\begin{document}

\title{Infinitely many two-variable generalisations\\of the
Alexander-Conway polynomial}
\asciititle{Infinitely many two-variable generalisations\\of the
Alexander-Conway polynomial}
\shorttitle{Two-variable generalisations of the
Alexander-Conway polynomial}

\authors{David De Wit\\Atsushi Ishii\\Jon Links}
\shortauthors{David De Wit, Atsushi Ishii and Jon Links}
\asciiauthors{David De Wit, Atsushi Ishii and Jon Links}

\addresses{{\rm DDW and JL:\qua}Department of Mathematics, The University of 
Queensland\\4072, Brisbane, Australia\\{\rm and}\\{\rm AI:\qua}Department 
of Mathematics, Graduate School of
Science, Osaka University\\Machikaneyama 1-16, Toyonaka, Osaka,
560-0043, Japan}

\gtemail{\href{mailto:Dr_David_De_Wit@yahoo.com.au}{Dr\_David\_De\_Wit@yahoo.com.au}, 
\mailto{aishii@cr.math.sci.osaka-u.ac.jp}, \mailto{jrl@maths.uq.edu.au}}

\asciiemail{Dr David De Wit@yahoo.com.au, 
aishii@cr.math.sci.osaka-u.ac.jp, jrl@maths.uq.edu.au}

\asciiaddress{DDW and JL: Department of Mathematics, The University of 
Queensland\\4072, Brisbane, Australia\\and\\AI: Department of
Mathematics, Graduate School of Science, Osaka
University\\Machikaneyama 1-16, Toyonaka, Osaka, 560-0043, Japan}

\begin{abstract}
We show that the Alexander-Conway polynomial $\Delta$ is obtainable via
a particular one-variable reduction of each two-variable Links--Gould
invariant $LG^{m,1}$, where $m$ is a positive integer.
Thus there exist infinitely many two-variable generalisations of
$\Delta$.
This result is not obvious since in the reduction, the representation of
the braid group generator used to define $LG^{m,1}$ does not satisfy a
second-order characteristic identity unless $m=1$.
To demonstrate that the one-variable reduction of $LG^{m,1}$ satisfies
the defining skein relation of $\Delta$, we evaluate the kernel of a
quantum trace.
\end{abstract}

\asciiabstract{%
We show that the Alexander-Conway polynomial Delta is obtainable via a
particular one-variable reduction of each two-variable Links-Gould
invariant LG^{m,1}, where m is a positive integer.  Thus there exist
infinitely many two-variable generalisations of Delta.  This result is
not obvious since in the reduction, the representation of the braid
group generator used to define LG^{m,1} does not satisfy a
second-order characteristic identity unless m=1.  To demonstrate that
the one-variable reduction of LG^{m,1} satisfies the defining skein
relation of Delta, we evaluate the kernel of a quantum trace.}


\asciikeywords{Link, knot, Alexander-Conway polynomial, Lie superalgebra,
Links-Gould invariant}

\maketitle

\section{Introduction}

The type I Lie superalgebras $sl(m|n)$ and $osp(2|2n)$ have the
distinguishing property that they admit nontrivial one-parameter
families of representations, and these representations extend to their
quantum deformations $U_q[sl(m|n)]$ and $U_q[osp(2|2n)]$.
Consequently, the link invariants derived from such representations are
two-variable invariants~\cite{GouldLinksZhang96,LinksGould92}.
In the simplest case $sl(1|1)$, the invariant reduces to a one-variable
invariant which is precisely the Alexander-Conway polynomial
$\Delta$~\cite{KauffmanSaleur91}.
The simplest nontrivial example of a two-variable invariant comes from
$sl(2|1)\cong{}osp(2|2)$~\cite{DeWitKauffmanLinks99,Ishii04,Ishii03skein,IshiiKanenobu03,LinksGould92}.
For this case it has recently been shown~\cite{Ishii03LGAC} that a
certain one-variable reduction recovers $\Delta$.
Whilst it may appear that the origin of this result may lie in the quantum
superalgebra embedding $U_q[sl(1|1)] \subset U_q[sl(2|1)]$, in fact
$\Delta$ is recovered only when the variable $q$ assumes specific roots
of unity.
It is also well known that $\Delta$ occurs as a one-variable reduction
of the two-variable HOMFLY polynomial~\cite{HOMFLY}.
The result of~\cite{Ishii03LGAC} thus shows that the extension of
$\Delta$ to a two-variable quantum invariant is not unique.

In this paper we extend the result of~\cite{Ishii03LGAC} to higher rank
superalgebras.
Specifically we employ $U_q[gl(m|1)]$, which differs from $U_q[sl(m|1)]$
by the addition of a central element; that is
$U_q[gl(m|1)] =U_q[u(1)\oplus sl(m|1)]$.
This yields the same link invariant, but conveniently makes the
representation theory easier to handle.
For the minimal one-parameter family of representations of dimension
$2^m$, we construct a link invariant denoted $LG^{m,1}(\tau,q)$ which is
a function of two independent variables $q$ and
$\tau\equiv q^{-\alpha}$.
Here, $\alpha$ is the complex parameter which indexes the underlying
representations.
These invariants have been introduced and studied
in~\cite{DeWit01,DeWit02,GouldLinksZhang96}.

Our main result is Theorem~\ref{thm:LGm1AC} (originally conjectured
in~\cite{Ishii03LGAC}), which is the following relation between
$LG^{m,1}$ and $\Delta$.
For an oriented link $L$, we have:
\begin{eqnarray}
LG^{m,1}_L(\tau,e^{\pi\sqrt{-1}/m})
= \Delta_L(\tau^{2m}).
\label{eq:LGm1AC}
\end{eqnarray}
We prove this relation by showing that
$LG^{m,1}(\tau,e^{\pi\sqrt{-1}/m})$ satisfies the skein relation
defining $\Delta(\tau^{2m})$.
To that end, we begin by recalling the method of construction for
$LG^{m,1}$, following~\cite{DeWit01,DeWit02}.
Next, we demonstrate a couple of technical lemmas from the
representation theory of $U_q[gl(m|1)]$.
Using them, the key to the proof involves determining the kernel of a
quantum trace, as per the method in~\cite{Ishii03LGAC}.
We stress that the representation of the braid group generator used in
the definition of $LG^{m,1}(\tau,e^{\pi\sqrt{-1}/m})$ does not satisfy a
second-order characteristic identity.
If this were the case, a proof would be trivial.
We also stress that, as for the $m=2$ case, our result does not directly
arise from the quantum superalgebra embedding
$U_q[sl(1|1)] \subset U_q[sl(m|1)]$.

\section{Quantum link invariants and $LG^{m,n}$}

Any oriented tangle diagram can be expressed up to isotopy as a diagram
composed from copies of the following elementary oriented tangle
diagrams.
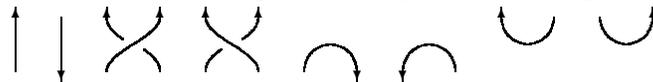
\begin{center}
\begin{picture}(0,20)
  \qbezier(0,0)(0,0)(0,20)
  \put(0.25,20){\vector(0,1){4}}
\end{picture}
\hspace{10pt}
\begin{picture}(0,20)
  \qbezier(0,0)(0,0)(0,20)
  \put(0.25,0){\vector(0,-1){4}}
\end{picture}
\hspace{10pt}
\begin{picture}(20,20)
  \qbezier(0,0)(0,4)(10,10)
  \qbezier(10,10)(20,16)(20,20)
  \qbezier(0,20)(0,16)(6,12)
  \qbezier(14,8)(20,4)(20,0)
  \put(0.25,20){\vector(0,1){4}}
  \put(20.25,20){\vector(0,1){4}}
\end{picture}
\hspace{10pt}
\begin{picture}(20,20)
  \qbezier(20,0)(20,4)(10,10)
  \qbezier(10,10)(0,16)(0,20)
  \qbezier(20,20)(20,16)(14,12)
  \qbezier(6,8)(0,4)(0,0)
  \put(0.25,20){\vector(0,1){4}}
  \put(20.25,20){\vector(0,1){4}}
\end{picture}
\hspace{10pt}
\begin{picture}(20,20)
  \qbezier(0,0)(1,9)(10,10)
  \qbezier(10,10)(19,9)(20,0)
  \put(20.25,0){\vector(0,-1){4}}
\end{picture}
\hspace{10pt}
\begin{picture}(20,20)
  \qbezier(0,0)(1,9)(10,10)
  \qbezier(10,10)(19,9)(20,0)
  \put(0.25,0){\vector(0,-1){4}}
\end{picture}
\hspace{10pt}
\begin{picture}(20,20)
  \qbezier(0,20)(1,11)(10,10)
  \qbezier(10,10)(19,11)(20,20)
  \put(0.25,20){\vector(0,1){4}}
\end{picture}
\hspace{10pt}
\begin{picture}(20,20)
  \qbezier(0,20)(1,11)(10,10)
  \qbezier(10,10)(19,11)(20,20)
  \put(20.25,20){\vector(0,1){4}}
\end{picture}
\end{center}
Furthermore any oriented tangle diagram can be expressed up to isotopy
as a \emph{sliced diagram} which is such a diagram sliced by horizontal
lines such that each domain between adjacent horizontal lines contains
either a single crossing or a single critical point.

Now let $V$ be a finite-dimensional vector space, with dual space
$V^*$.
Using these, we assign an invertible endomorphism
$R: V \otimes V \to V \otimes V$ and linear maps
$n: V \otimes V^* \to \mathbb{C}$,
$\tilde{n}: V^* \otimes V \to \mathbb{C}$,
$u: \mathbb{C} \to V^* \otimes V$
and
$\tilde{u}: \mathbb{C} \to V \otimes V^*$
to the elementary oriented tangle diagrams, as follows.
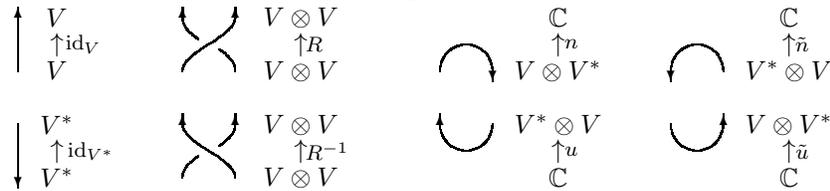
\begin{center}
\begin{picture}(40,20)\small
  \qbezier(0,0)(0,0)(0,20)
  \put(0.25,20){\vector(0,1){4}}
  \put(15,0){\makebox(0,0){$V$}}
  \put(15,10){\makebox(0,0){$\uparrow$}}
  \put(25,10){\makebox(0,0){\scriptsize$\mathrm{id}_V$}}
  \put(15,20){\makebox(0,0){$V$}}
\end{picture}
\hspace{15pt}
\begin{picture}(75,20)\small
  \qbezier(0,0)(0,4)(10,10)
  \qbezier(10,10)(20,16)(20,20)
  \qbezier(0,20)(0,16)(6,12)
  \qbezier(14,8)(20,4)(20,0)
  \put(0.25,20){\vector(0,1){4}}
  \put(20.25,20){\vector(0,1){4}}
  \put(45,0){\makebox(0,0){$V \otimes V$}}
  \put(45,10){\makebox(0,0){$\uparrow$}}
  \put(50,10){\makebox(0,0){\scriptsize$R$}}
  \put(45,20){\makebox(0,0){$V \otimes V$}}
\end{picture}
\hspace{15pt}
\begin{picture}(65,20)\small
  \qbezier(0,0)(1,9)(10,10)
  \qbezier(10,10)(19,9)(20,0)
  \put(20.25,0){\vector(0,-1){4}}
  \put(45,20){\makebox(0,0){$\mathbb{C}$}}
  \put(45,10){\makebox(0,0){$\uparrow$}}
  \put(50,10){\makebox(0,0){\scriptsize$n$}}
  \put(45,0){\makebox(0,0){$V \otimes V^*$}}
\end{picture}
\hspace{15pt}
\begin{picture}(65,20)\small
  \qbezier(0,0)(1,9)(10,10)
  \qbezier(10,10)(19,9)(20,0)
  \put(0.25,0){\vector(0,-1){4}}
  \put(45,20){\makebox(0,0){$\mathbb{C}$}}
  \put(45,10){\makebox(0,0){$\uparrow$}}
  \put(50,10){\makebox(0,0){\scriptsize$\tilde{n}$}}
  \put(45,0){\makebox(0,0){$V^* \otimes V$}}
\end{picture}
\end{center}
\vspace{10pt}
\begin{center}
\begin{picture}(40,20)\small
  \qbezier(0,0)(0,0)(0,20)
  \put(0.25,0){\vector(0,-1){4}}
  \put(15,0){\makebox(0,0){$V^*$}}
  \put(15,10){\makebox(0,0){$\uparrow$}}
  \put(28,10){\makebox(0,0){\scriptsize$\mathrm{id}_{V^*}$}}
  \put(15,20){\makebox(0,0){$V^*$}}
\end{picture}
\hspace{15pt}
\begin{picture}(75,20)\small
  \qbezier(20,0)(20,4)(10,10)
  \qbezier(10,10)(0,16)(0,20)
  \qbezier(20,20)(20,16)(14,12)
  \qbezier(6,8)(0,4)(0,0)
  \put(0.25,20){\vector(0,1){4}}
  \put(20.25,20){\vector(0,1){4}}
  \put(45,0){\makebox(0,0){$V \otimes V$}}
  \put(45,10){\makebox(0,0){$\uparrow$}}
  \put(55,10){\makebox(0,0){\scriptsize$R^{-1}$}}
  \put(45,20){\makebox(0,0){$V \otimes V$}}
\end{picture}
\hspace{15pt}
\begin{picture}(65,20)\small
  \qbezier(0,20)(1,11)(10,10)
  \qbezier(10,10)(19,11)(20,20)
  \put(0.25,20){\vector(0,1){4}}
  \put(45,0){\makebox(0,0){$\mathbb{C}$}}
  \put(45,10){\makebox(0,0){$\uparrow$}}
  \put(50,10){\makebox(0,0){\scriptsize$u$}}
  \put(45,20){\makebox(0,0){$V^* \otimes V$}}
\end{picture}
\hspace{15pt}
\begin{picture}(65,20)\small
  \qbezier(0,20)(1,11)(10,10)
  \qbezier(10,10)(19,11)(20,20)
  \put(20.25,20){\vector(0,1){4}}
  \put(45,0){\makebox(0,0){$\mathbb{C}$}}
  \put(45,10){\makebox(0,0){$\uparrow$}}
  \put(50,10){\makebox(0,0){\scriptsize$\tilde{u}$}}
  \put(45,20){\makebox(0,0){$V \otimes V^*$}}
\end{picture}
\end{center}
Corresponding to an oriented tangle diagram $D$, we then obtain a linear
map $[D]$ by composing tensor products of copies of the linear maps
associated with the elementary tangle diagrams in $D$. For example:
\begin{eqnarray}
\left[ \hspace{1.5pt}
\begin{minipage}{20pt}
\begin{picture}(20,34)
  \qbezier(0,10)(0,12)(5,15)
  \qbezier(5,15)(10,18)(10,20)
  \qbezier(0,20)(0,18)(3,16)
  \qbezier(7,14)(10,12)(10,10)
  \qbezier(10,20)(10.5,24.5)(15,25)
  \qbezier(15,25)(19.5,24.5)(20,20)
  \qbezier(10,10)(10.5,5.5)(15,5)
  \qbezier(15,5)(19.5,5.5)(20,10)
  \qbezier(0,10)(0,0)(0,0)
  \qbezier(0,30)(0,20)(0,20)
  \qbezier(20,20)(20,10)(20,10)
  \put(.25,30){\vector(0,1){4}}
\end{picture}
\end{minipage}
\hspace{1.5pt} \right]
=(\mathrm{id}_V \otimes n)
 (R \otimes \mathrm{id}_{V^*})
 (\mathrm{id}_V \otimes u).
\label{eq:bracketofR}
\end{eqnarray}
A quantum link invariant may then be defined as follows.
Set $V$ as the module associated with an irreducible, finite-dimensional
representation $\pi$ of some ribbon Hopf (super)algebra, for instance a
quantum superalgebra.
We then obtain the bracket $[~]$ by setting $R$ as a representation of
the braid group generator associated with the tensor product
representation $\pi\otimes\pi$.
This choice ensures the invariance of the bracket under the second
Reidemeister move, due to the invertibility of $R$, and the third
Reidemeister move, as $R$ satisfies the Yang--Baxter equation
(see~(\ref{eq:yangbaxterequation}) below).
Note that at this point, we may freely use any scaling of $R$.

Now let the \emph{quantum trace} be the linear map
$\mathrm{cl}:
\mathrm{End}(V^{\otimes (k+1)}) \to \mathrm{End}(V^{\otimes k})$
(where $k \geqslant 1$), which is defined for
$X\in\mathrm{End}(V^{\otimes (k+1)})$ by:
\begin{eqnarray*}
\mathrm{cl}(X)
=(\mathrm{id}_V^{\otimes k} \otimes n)
 (X \otimes \mathrm{id}_{V^*})
 (\mathrm{id}_V^{\otimes k} \otimes u).
\end{eqnarray*}
Observe that~(\ref{eq:bracketofR}) describes $\mathrm{cl}(R)$.
Demanding that $\mathrm{cl}(R) =\mathrm{cl}(R^{-1}) =\mathrm{id}_V$
ensures the invariance of the bracket under the first Reidemeister
move.
This requirement determines the scaling of $R$, and also the choice of
the mappings $n$, $\tilde{n}$, $u$ and $\tilde{u}$.
Specifically, representation-theoretic considerations mean that these
mappings may be defined in terms of the representation of an element of
the Cartan subalgebra of the underlying (super)algebra
(see~\cite{DeWit01}).

For any given oriented tangle $T$, we thus obtain a map $[D_T]$, where
$D_T$ is an oriented tangle diagram corresponding to $T$, and the map
$[D_T]$ is invariant under ambient isotopy of $T$.
For notational convenience, we shall generally write $[T]$ for $[D_T]$,
and this is meaningful as the evaluation of the invariant is independent
of the choice of diagram $D_T$.
By construction, the maps $R$, $R^{-1}$, $n$, $\tilde{n}$, $u$ and
$\tilde{u}$ are invariant with respect to the action of the Hopf
(super)algebra.
Consequently, the map $[T]$ is also invariant with respect to this
action.
Specifically, where $T$ is an oriented $(1,1)$-tangle, the choice of $V$
as irreducible means that Schur's Lemma ensures that $[T]$ is a scalar
map (that is, a scalar multiple of $\mathrm{id}_V$).
This scalar is then a quantum link invariant of the link $\widehat{T}$
formed by the closure of $T$ (see~\cite{Ohtsuki02,Zhang95}); in
particular the scalar is unity when $\widehat{T}$ is the unknot.

\vspace{\baselineskip}

Now fix positive integers $m$ and $n$, and consider the quantum
superalgebra $U_q[gl(m|n)]$, a quantum deformation of the universal
enveloping algebra of the Lie superalgebra $gl(m|n)$.
The two-variable Links--Gould invariant $LG^{m,n}(\tau,q)$ may then be
obtained by specialising the above framework to the case of the minimal
$2^{mn}$-dimensional $U_q[gl(m|n)]$ representation $\pi$ bearing a free
parameter $\alpha$ (for details,
see~\cite{DeWit01,DeWitKauffmanLinks99}).
In that case, where $V$ is the module associated with $\pi$, we
explicitly write:
\begin{eqnarray}
[T] =LG^{m,n}_{\widehat{T}}(\tau,q) \, \mathrm{id}_V,
\label{eq:LGdefinition}
\end{eqnarray}
where we have used the variable $\tau=q^{-\alpha}$ instead of $\alpha$;
below we freely interchange use of the variables $\alpha$ and $\tau$.
Note that we have $LG^{m,n}_{\bigcirc}(\tau,q)=1$.

Next, we present an important symmetry of these invariants.
To that end, firstly note that $U_q[gl(m|n)]$ is defined
(see~\cite{Zhang93}) in terms of a fixed invariant bilinear form on the
weight space of $gl(m|n)$.
We adopt the convention that the form is positive definite for $gl(m)$
roots and negative definite for $gl(n)$ roots.
It may be deduced from the definition that, under this convention, the
following superalgebra isomorphism holds:
\begin{eqnarray}
U_q[gl(m|n)] \cong U_{q^{-1}}[gl(n|m)].
\label{eq:UqglmncongUqglnm}
\end{eqnarray}
We then note that the substitution $\alpha \rightarrow -(\alpha+m-n)$
maps the $U_q[gl(m|n)]$ representation $\pi$ to its dual $\pi^*$.
This, together with~(\ref{eq:UqglmncongUqglnm}) allows us to deduce
that, for any oriented link $L$, we have:
\begin{eqnarray}
LG_L^{m,n}(\tau,q) = LG_L^{n,m}(\tau,q^{-1}).
\label{eq:SymmetryofLGmn}
\end{eqnarray}
We shall be interested below in the case $LG^{m,1}$ and the substitution
of the root of unity $e^{\pi \sqrt{-1}/m}$ for $q$; importantly, the
structure of the representation does not change at this particular root
of unity.
We also emphasise that under this substitution, we intend $\tau$ to
remain independent; that is, we do \emph{not} express it as
$e^{-\alpha\pi\sqrt{-1}/m}$.

\section{Some $U_q[gl(m|1)]$ representation theory}

The construction of the mappings $R$, $R^{-1}$, $n$, $\tilde{n}$, $u$
and $\tilde{u}$ determining $LG^{m,1}$ can be described in terms of the
representation theory of $U_q[gl(m|1)]$.
In this section, we establish notational conventions and provide the
necessary representation-theoretic results needed to deduce our main
result, relation~(\ref{eq:LGm1AC}).

We begin with the fact that every irreducible finite-dimensional
$U_q[gl(m|1)]$ module $V(\Lambda)$ is uniquely labelled by its highest
weight $\Lambda=(\Lambda_1,\dots,\Lambda_m|\Lambda_{m+1})$.
Moreover, each $V(\Lambda)$ is completely reducible with respect to the
even subalgebra $U_q[gl(m)\oplus gl(1)]$ such that we may write
\begin{eqnarray*}
V(\Lambda) = \bigoplus_k V^0(\mu_k),
\end{eqnarray*}
where each $V^0(\mu_k)$ is an irreducible $U_q[gl(m)\oplus gl(1)]$
module with highest weight $\mu_k$.
Here, we are in fact only interested in a subclass of these
$U_q[gl(m|1)]$ modules, that is those whose highest weights are of the
form
\begin{eqnarray*}
\Lambda(i,j,\alpha)
\triangleq ({0}_{m-i-j},{-1}_i,{-2}_j\,|\,\alpha+i+2j),
\end{eqnarray*}
where the subscripts indicate the number of times each entry is repeated
in the weight, and $\alpha$ is an arbitrary complex parameter.
We set $V(i,j,\alpha)$ as the irreducible module with highest weight
$\Lambda(i,j,\alpha)$, and we also let $V^0(i,j,\alpha)$ denote the
irreducible $U_q[gl(m)\oplus gl(1)]$ module with the same highest
weight.

Specifically, $LG^{m,1}$ is defined in terms of the representation
associated with the module $V(0,0,\alpha)$.
We have the following decompositions~\cite{GouldLinksZhang96}:
\begin{eqnarray}
V(0,0,\alpha)
&=& \bigoplus_{i=0}^m V^0(i,0,\alpha),
\label{eq:decompofV} \\
V(0,0,\alpha)\otimes V(0,0,\alpha)
&=& \bigoplus_{i=0}^m V(i,0,2\alpha).
\label{eq:decompofV00alphaV00alpha}
\end{eqnarray}
As each submodule $V(i,0,2\alpha)$
in~(\ref{eq:decompofV00alphaV00alpha}) is \emph{typical}, applying the
Kac induced module construction~\cite{Kac78}, we may similarly deduce
the following decomposition:
\begin{eqnarray}
V(i,0,2\alpha)
=\bigoplus_{j=0}^i \bigoplus_{k=i}^m V^0(k-j,j,2\alpha).
\label{eq:decompofVi02alpha}
\end{eqnarray}
In~\cite{GouldLinksZhang96} the decompositions~(\ref{eq:decompofV})
and~(\ref{eq:decompofV00alphaV00alpha}) were deduced for generic values
of $\alpha$ and real, positive $q$.
It is important to stress
that~(\ref{eq:decompofV})--(\ref{eq:decompofVi02alpha}) remain valid
when $q=e^{\pi\sqrt{-1}/m}$.
We comment further on this aspect in the proof of
Lemma~\ref{lem:welldefined} (below).

To simplify notation, we shall write $V$ for $V(0,0,\alpha)$.
With respect to~(\ref{eq:decompofV00alphaV00alpha}), setting $V_i$ as
$V(i,0,2\alpha)$, let $P_i$ be the projector mapping $V \otimes V$ onto
$V_i$, so that we have:
\begin{eqnarray}
P_i P_j =\delta_{ij} P_i,
\qquad \qquad
P_0 + \cdots + P_m =\mathrm{id}_{V\otimes V}.
\label{eq:id}
\end{eqnarray}
Then, from~\cite{DeWit02}, we have:
\begin{eqnarray}
R =\sum_{i=0}^m \xi_i P_i,
\qquad \qquad
R^{-1} =\sum_{i=0}^m \xi_i^{-1} P_i,
\label{eq:R}
\end{eqnarray}
where
\begin{eqnarray}
\xi_i
=(-1)^i q^{i(2\alpha +i-1) -m\alpha}
\equiv (-1)^i \tau^{m-2i} q^{i(i-1)}.
\label{eq:xi}
\end{eqnarray}
Note that the scaling of $R$ has been chosen such that
$\mathrm{cl}(R) =\mathrm{cl}(R^{-1}) =\mathrm{id}_V$.
The grading of the underlying vector space $V$ means that $R$ as defined
in~(\ref{eq:R}) actually satisfies a \emph{graded} Yang--Baxter
equation~\cite{Zhang95}.
However, by insertion of factors of $-1$ into some of the components of
$R$ (as described in~\cite{DeWit02}) it is made to satisfy the usual
ungraded Yang--Baxter equation:
\begin{eqnarray}
(R\otimes\mathrm{id}_V)(\mathrm{id}_V \otimes R)
(R\otimes\mathrm{id}_V)=(\mathrm{id}_V \otimes R)
(R\otimes\mathrm{id}_V)(\mathrm{id}_V \otimes R).
\label{eq:yangbaxterequation}
\end{eqnarray}
It is clear from~(\ref{eq:R}) that $R$ satisfies the characteristic
identity of order $m+1$:
\begin{eqnarray}
\prod_{i=0}^m (R -\xi_i \, \mathrm{id}_{V\otimes V})
= 0 \, \mathrm{id}_{V\otimes V}.
\label{eq:characteristicidentityofRm1}
\end{eqnarray}
For any linear map $X$ we denote $X|_{q=e^{\pi\sqrt{-1}/m}}$ by
$\overline{X}$.
Similarly, for any vector space $W$ over
$\mathbb{C}[q,q^{-1},\tau,\tau^{-1}]$, we denote by $\overline{W}$ the
vector space over $\mathbb{C}[\tau,\tau^{-1}]$ obtained from $W$ by
setting $q=e^{\pi\sqrt{-1}/m}$.
It is necessary to affirm that the mappings $R$, $R^{-1}$, $n$,
$\tilde{n}$, $u$, $\tilde{u}$ and each $P_i$ are well-defined in the
substitution $q=e^{\pi\sqrt{-1}/m}$.

\begin{lemma} \label{lem:welldefined}
The mappings $\overline{R}$, $\overline{R^{-1}}$, $\overline{n}$,
$\overline{\tilde{n}}$, $\overline{u}$, $\overline{\tilde{u}}$ and each
$\overline{P_i}$ are well-defined, that is, all matrix elements of $R$,
$R^{-1}$, $n$, $\tilde{n}$, $u$, $\tilde{u}$ and each $P_i$ have no pole
at $q=e^{\pi\sqrt{-1}/m}$.
\end{lemma}

\begin{proof}
We begin by recalling from~\S3 of~\cite{LinksGouldZhang93} the
$U_q[gl(m|1)]$ central element
$\Gamma\triangleq(v\otimes v)\Delta(v^{-1})$, where $\Delta$ is the
coproduct and $v$ is the ribbon element in the centre of
$U_q[gl(m|1)]$.
Each projector $P_i$ may be expressed as a polynomial function of the
representation of $\Gamma$ via:
\begin{eqnarray*}
P_i =\prod_{j\neq i}
\frac{\Gamma -\gamma_j \, \mathrm{id}_{V\otimes V}}{\gamma_i -\gamma_j},
\end{eqnarray*}
where $\gamma_i$ denotes the eigenvalue of $\Gamma$ on $V_i$.
In fact, $\gamma_i=\xi_i^2$, where $\xi_i$ is as introduced
in~(\ref{eq:xi}).
Note that $\overline{\gamma_i}\neq\overline{\gamma_j}$, for $i\neq j$.
If we rewrite $P_i =N_i/D_i$, where:
\begin{eqnarray*}
N_i \triangleq
\prod_{j\neq i} \left( \Gamma -\gamma_j \, \mathrm{id}_{V\otimes V} \right),
\qquad \mathrm{and} \qquad
D_i \triangleq
\prod_{j\neq i} \left( \gamma_i -\gamma_j \right),
\end{eqnarray*}
then we see that $\overline{D_i}$ is nonzero.
We next show that $\overline{N_i}$ is well-defined.

To that end, $\Gamma$ may be expressed as a power series over
$\mathbb{C}[q,q^{-1}]$ of the simple $U_q[gl(m|1)]$
generators~\cite[\S3]{LinksGouldZhang93}.
As the weights $\Lambda(i,j,\alpha)$ are generically
\emph{essentially typical}, for general $q$, expressions are
known~\cite{PalevStoilovaVanderJeugt94} for the matrix elements of the
simple generators in a Gel'fand--Zetlin basis.
The matrix elements of the even simple generators are well-defined when
$q=e^{\pi\sqrt{-1}/m}$.
This follows as condition~(3.2)
of~\cite{AbdesselamArnaudonChakrabarti95} is satisfied for all the
modules $V^0(i,0,\alpha)$ of (\ref{eq:decompofV}).
Thus, (\ref{eq:decompofV}), and by the same reasoning
(\ref{eq:decompofVi02alpha}), remains valid for $q=e^{\pi\sqrt{-1}/m}$.
The matrix elements of the odd simple generators for $U_q[gl(m|1)]$ are
given by formulae (27,28) of~\cite{PalevStoilovaVanderJeugt94}.
Unlike the situation for the even generators, these formulae explicitly
depend on the variable $\alpha$.
This means that they are well-defined when $q=e^{\pi\sqrt{-1}/m}$, since
their denominators are nonvanishing for generic values of $\alpha$.

Thus, each $\overline{N_i}$, hence each $\overline{P_i}$ is
well-defined, and consequently so are $\overline{R}$ and
$\overline{R^{-1}}$.
The fact that the mappings $\overline{n}$, $\overline{\tilde{n}}$,
$\overline{u}$ and $\overline{\tilde{u}}$ are also well-defined follows
from their definitions in terms of the representation of an element of
the $U_q[gl(m|1)]$ Cartan subalgebra.
\end{proof}

We remark that this proof also demonstrates that the decomposition
of~(\ref{eq:decompofV00alphaV00alpha}) remains valid in the reduction
$q=e^{\pi\sqrt{-1}/m}$, since the projectors remain well-defined.

\begin{lemma} \label{lem:P_i}
For each $i=0,\dots, m$ the expression $\overline{\mathrm{cl}(P_i)}$ is
a well-defined scalar multiple of $\mathrm{id}_{\overline{V}}$; in fact
$\overline{\mathrm{cl}(P_i)}=0 \, \mathrm{id}_{\overline{V}}$ for
$i=1,\dots,m-1$.
\end{lemma}

\begin{proof}
Theorem 1 of~\cite{GouldLinksZhang96}, specified to our situation,
reads:
\begin{eqnarray*}
&& \mathrm{cl}(P_i) =(-1)^i \prod_{j=1}^i
\frac{q^{m-j+1} -q^{-(m-j+1)}}{q^{i-j+1} -q^{-(i-j+1)}}
\cdot
\frac{\tau q^{-(j-1)} -\tau^{-1} q^{j-1}}
     {\tau^2 q^{-(i+j-2)} -\tau^{-2} q^{i+j-2}}
\nonumber \\
&& \hspace{100pt}
\times \prod_{j=i+1}^m
\frac{\tau q^{-(j-1)} -\tau^{-1} q^{j-1}}
     {\tau^2 q^{-(i+j-1)} -\tau^{-2} q^{i+j-1}} \, \mathrm{id}_V;
\end{eqnarray*}
note that in the cases $i=0,m$, the formula reduces to the following:
\begin{eqnarray*}
\mathrm{cl}(P_0) &=& \prod_{j=1}^m
\frac{\tau q^{-(j-1)} -\tau^{-1} q^{j-1}}
     {\tau^2 q^{-(j-1)} -\tau^{-2} q^{j-1}} \, \mathrm{id}_V \\
\mathrm{cl}(P_m) &=& (-1)^m \prod_{j=1}^m
\frac{\tau q^{-(j-1)} -\tau^{-1} q^{j-1}}
     {\tau^{2} q^{-(m+j-2)} -\tau^{-2} q^{m+j-2}} \, \mathrm{id}_V.
\end{eqnarray*}
In these formulae, we intend $\tau\equiv q^{-\alpha}$ to be restricted
so that the complex variable $\alpha$ is not an integer.
(By an analytic continuation argument, this restriction does not affect
our final result.)
Now observe that the denominator of $\mathrm{cl}(P_i)$ never contains
any factors of $q^m-q^{-m}$; this means that
$\overline{\mathrm{cl}(P_i)}$ is always well-defined.
However, if $i\neq 0,m$, its numerator always contains a factor of
$q^m-q^{-m}$, meaning that
$\overline{\mathrm{cl}(P_i)}=0 \, \mathrm{id}_{\overline{V}}$.
\end{proof}

Now let $V$ have a weight basis $\{e_0,\dots,e_{2^m-1}\}$.
Since the weight spectrum of $V$ is multiplicity-free, we can choose the
labelling such that for $i=0,\dots,m$, the vector $e_i$ has weight
$\Lambda(i,0,\alpha)$.
In terms of this basis, any $A\in\mathrm{End}(V\otimes V)$ may be
written in component form via
$A (e_k\otimes e_l) =\sum_{ij} A^{ij}_{kl} (e_i\otimes e_j)$.

\begin{lemma} \label{lem:Pijjjjsdeltaij}
$(P_i)^{jj}_{jj}=\delta_{ij}$, for all $i,j=0,\dots,m$.
\end{lemma}

\begin{proof}
From~(\ref{eq:decompofV}), we know that $e_i$ is a
$U_q[gl(m)\oplus gl(1)]$ highest weight vector.
Therefore $v_i\triangleq e_i\otimes e_i$ is also a
$U_q[gl(m)\oplus gl(1)]$ highest weight vector, of weight
$\Lambda (0,i,2\alpha)$.
Now looking at~(\ref{eq:decompofVi02alpha}), we see that this
$U_q[gl(m)\oplus gl(1)]$ highest weight only occurs in $V(i,0,2\alpha)$.
Thus, $v_i$ generates the irreducible module $V(i,0,2\alpha)$, and
moreover, for each $V(j,0,2\alpha)$ there exists a
$v_j\equiv e_j\otimes e_j$ which generates it.
Thus, for each projector $P_i$ we have
$P_i(e_j\otimes e_j) =\delta_{ij} (e_j\otimes e_j)$,
hence we conclude $(P_i)^{jj}_{jj} =\delta_{ij}$.
\end{proof}

\section{The relation}

In this section we show the following relation:
\begin{eqnarray*}
LG^{m,1}_L(\tau,e^{\pi\sqrt{-1}/m}) =\Delta_L(\tau^{2m}),
\end{eqnarray*}
where $\Delta_L(t)$ is the Alexander-Conway polynomial which is defined
by the following relations:
\begin{eqnarray}
\Delta_{\bigcirc}(t) &=& 1,
\label{eq:ACunknotisunity} \\
\Delta_{
\begin{minipage}{10pt}
\begin{picture}(10,10)
  \qbezier(0,0)(0,0)(10,10)
  \qbezier(10,0)(10,0)(6.5,3.5)
  \qbezier(0,10)(0,10)(3.5,6.5)
  \put(-.8,11){\vector(-1,1){0}}
  \put(11.2,11){\vector(1,1){0}}
\end{picture}
\end{minipage}
}(t) -\Delta_{
\begin{minipage}{10pt}
\begin{picture}(10,10)
  \qbezier(0,10)(0,10)(10,0)
  \qbezier(0,0)(0,0)(3.5,3.5)
  \qbezier(6.5,6.5)(10,10)(10,10)
  \put(-.8,11){\vector(-1,1){0}}
  \put(11.2,11){\vector(1,1){0}}
\end{picture}
\end{minipage}
}(t) &=& (t^{1/2}-t^{-1/2}) \Delta_{
\begin{minipage}{10pt}
\begin{picture}(10,10)
  \qbezier(0,0)(6,5)(0,10)
  \qbezier(10,0)(4,5)(10,10)
  \put(-.8,11){\vector(-1,1){0}}
  \put(11.2,11){\vector(1,1){0}}
\end{picture}
\end{minipage}
}(t). \label{eq:ACskein}
\end{eqnarray}

\begin{lemma} \label{lem:T}
Where $T$ is an oriented $(2,2)$-tangle, $[T]$ may be expressed as:
\begin{eqnarray}
[T] =\sum_{i=0}^m a_i^T P_i,
\label{eq:TisaiPi}
\end{eqnarray}
where the coefficients $a_i^T$ are such that each $\overline{a_i^T}$ is
well-defined.
\end{lemma}

\begin{proof}
Firstly, note that $[T]$ is a product of $U_q[gl(m|1)]$-invariant
mappings, and $\{P_0,\dots,P_m\}$ is a basis for the space of such
mappings on $V\otimes V$.
Thus, $[T]$ is necessarily of the form~(\ref{eq:TisaiPi}).
Recall from Lemma~\ref{lem:welldefined} that the mappings $R$, $R^{-1}$,
$n$, $\tilde{n}$, $u$ and $\tilde{u}$ are well-defined in the
substitution $q=e^{\pi\sqrt{-1}/m}$.
Thus, as $[T]$ is defined in terms of these mappings, it is also
well-defined in the substitution.
Then, using Lemma~\ref{lem:Pijjjjsdeltaij}, we have:
\begin{eqnarray*}
[T]^{jj}_{jj} =\sum_{i=0}^m a^T_i (P_i)^{jj}_{jj} =a^T_j,
\qquad \mathrm{for} \; j=0,\dots,m,
\end{eqnarray*}
and conversely, for each index $i=0,\dots,m$, we have
$a^T_i=[T]^{ii}_{ii}$.
Thus, as $\overline{[T]}$ is well-defined, so is $\overline{a_i^T}$.
\end{proof}

Before moving on to our main result, we emphasise that $\overline{R}$
does not satisfy a second-order characteristic identity (unless $m=1$).
In particular, the following identity:
\begin{eqnarray}
\overline{R} -\overline{R^{-1}} =
(\tau^{m} -\tau^{-m}) \, \mathrm{id}_{\overline{V}\otimes\overline{V}},
\label{eq:2ndorderidentity}
\end{eqnarray}
\emph{only} holds for $m=1$.
If~(\ref{eq:2ndorderidentity}) held for arbitrary $m$, the proof of our
main result would be trivial.

\begin{theorem} \label{thm:LGm1AC}
For any oriented link $L$, there holds:
\begin{eqnarray*}
LG^{m,1}_L(\tau,e^{\pi\sqrt{-1}/m}) =\Delta_L(\tau^{2m}).
\end{eqnarray*}
\end{theorem}

\begin{proof}
For any oriented $(2,2)$-tangle $T$, we have:
\def\strut{\vrule width0pt height 15pt}
\begin{eqnarray*}
\lefteqn{
\overline{ \left[ \hspace{5pt}
\begin{minipage}{30pt}
\begin{picture}(30,62)
  \qbezier(15,25)(20,30)(13,37)
  \qbezier(5,25)(0,30)(7,37)
  \qbezier(7,37)(10,40)(13,43)
  \put(10,20){\circle{14}}
  \put(10,20){\makebox(0,0){$T$}}
  \qbezier(5,15)(0,10)(0,0)
  \qbezier(15,15)(20,10)(20,10)
  \qbezier(20,10)(30,0)(30,20)
  \qbezier(30,20)(30,20)(30,40)
  \qbezier(30,40)(30,60)(20,50)
  \qbezier(20,50)(20,50)(13,43)
  \put(0.3,60){\vector(0,1){1}}
  \qbezier(7,43)(0,50)(0,60)
  \put(30.3,30){\vector(0,-1){1}}
\end{picture}
\end{minipage}
\hspace{5pt} \right]}
-\overline{ \left[ \hspace{5pt}
\begin{minipage}{30pt}
\begin{picture}(30,62)
  \qbezier(5,25)(0,30)(7,37)
  \qbezier(15,25)(20,30)(13,37)
  \qbezier(13,37)(10,40)(7,43)
  \put(10,20){\circle{14}}
  \put(10,20){\makebox(0,0){$T$}}
  \qbezier(5,15)(0,10)(0,0)
  \qbezier(15,15)(20,10)(20,10)
  \qbezier(20,10)(30,0)(30,20)
  \qbezier(30,20)(30,20)(30,40)
  \qbezier(30,40)(30,60)(20,50)
  \qbezier(20,50)(20,50)(13,43)
  \put(0.3,60){\vector(0,1){1}}
  \qbezier(7,43)(0,50)(0,60)
  \put(30.3,30){\vector(0,-1){1}}
\end{picture}
\end{minipage}
\hspace{5pt} \right]}
-(\tau^m-\tau^{-m}) \overline{ \left[ \hspace{5pt}
\begin{minipage}{30pt}
\begin{picture}(30,62)
  \put(10,20){\circle{14}}
  \put(10,20){\makebox(0,0){$T$}}
  \qbezier(5,15)(0,10)(0,0)
  \qbezier(15,15)(20,10)(20,10)
  \qbezier(20,10)(30,0)(30,20)
  \qbezier(15,25)(20,30)(20,30)
  \qbezier(20,30)(30,40)(30,20)
  \put(0.3,60){\vector(0,1){1}}
  \qbezier(5,25)(0,30)(0,60)
  \put(30.3,20){\vector(0,-1){1}}
\end{picture}
\end{minipage}
\hspace{5pt} \right]}
} \\
&=&\strut
 \overline{\mathrm{cl}(R \circ [T])}
-\overline{\mathrm{cl}(R^{-1} \circ [T])}
-(\tau^m-\tau^{-m})\overline{\mathrm{cl}([T])} \\
&=&\strut
 \overline{{\textstyle\sum_{i=0}^m} \xi_i a_i^T \mathrm{cl}(P_i)}
-\overline{{\textstyle\sum_{i=0}^m} \xi_i^{-1} a_i^T \mathrm{cl}(P_i)}
-(\tau^m-\tau^{-m})
 \overline{{\textstyle\sum_{i=0}^m}a_i^T\mathrm{cl}(P_i)} \\
&=&\strut
 {\textstyle\sum_{i=0}^m}
  \overline{\xi_i} \ \overline{a_i^T} \ \overline{\mathrm{cl}(P_i)}
-{\textstyle\sum_{i=0}^m}
  \overline{\xi_i^{-1}} \ \overline{a_i^T} \ \overline{\mathrm{cl}(P_i)}
-(\tau^m-\tau^{-m}){\textstyle\sum_{i=0}^m}
  \overline{a_i^T} \ \overline{\mathrm{cl}(P_i)} \\
&=&\strut
 {\textstyle\sum_{i=0}^m} \left(
 (\overline{\xi_i} -\overline{\xi_i^{-1}}) -(\tau^m-\tau^{-m})
 \right)
 \overline{a_i^T} \ \overline{\mathrm{cl}(P_i)} \\
&=&\strut
0 \, \mathrm{id}_{\overline{V}},
\end{eqnarray*}
where the second equality follows from~(\ref{eq:id}), (\ref{eq:R})
and~(\ref{eq:TisaiPi}), the third from Lemmas~\ref{lem:P_i}
and~\ref{lem:T}, and the last from Lemma~\ref{lem:P_i} and the
observation from~(\ref{eq:xi}) that
\begin{eqnarray}
\overline{\xi_0} =-\overline{\xi_m^{-1}} =\tau^m.
\label{eq:xi0isximistaum}
\end{eqnarray}
In view of~(\ref{eq:LGdefinition}), we thus have the following skein
relation:
\begin{eqnarray*}
\overline{ LG^{m,1}_{ \hspace*{1pt}
\begin{minipage}{10pt}
\begin{picture}(10,10)
  \qbezier(0,0)(0,0)(10,10)
  \qbezier(10,0)(10,0)(6.5,3.5)
  \qbezier(0,10)(0,10)(3.5,6.5)
  \put(-.8,11){\vector(-1,1){0}}
  \put(11.2,11){\vector(1,1){0}}
\end{picture}
\end{minipage}
\hspace*{1.5pt} } }
-\overline{ LG^{m,1}_{ \hspace*{1pt}
\begin{minipage}{10pt}
\begin{picture}(10,10)
  \qbezier(0,10)(0,10)(10,0)
  \qbezier(0,0)(0,0)(3.5,3.5)
  \qbezier(6.5,6.5)(10,10)(10,10)
  \put(-.8,11){\vector(-1,1){0}}
  \put(11.2,11){\vector(1,1){0}}
\end{picture}
\end{minipage}
\hspace*{1.5pt} } }
 =(\tau^m-\tau^{-m})
\overline{ LG^{m,1}_{ \hspace*{1pt}
\begin{minipage}{10pt}
\begin{picture}(10,10)
  \qbezier(0,0)(6,5)(0,10)
  \qbezier(10,0)(4,5)(10,10)
  \put(-.8,11){\vector(-1,1){0}}
  \put(11.2,11){\vector(1,1){0}}
\end{picture}
\end{minipage}
\hspace*{1.5pt} } },
\end{eqnarray*}
so $\overline{LG^{m,1}}$ satisfies~(\ref{eq:ACskein}). It also
satisfies~(\ref{eq:ACunknotisunity}), as
$\overline{LG^{m,1}_{\bigcirc}}=1$.
Thus, for any oriented link $L$, we have
$LG^{m,1}_L(\tau,e^{\pi\sqrt{-1}/m}) =\Delta_L(\tau^{2m})$.
\end{proof}

Now note that the proof of Theorem~\ref{thm:LGm1AC} remains valid when
$\overline{X}$ is instead regarded as $X|_{q=e^{\pi\sqrt{-1}\,r/m}}$,
where $r$ is any integer such that $r$ and $m$ are relatively prime.
This follows since Lemmas~\ref{lem:welldefined}, \ref{lem:P_i}
and~\ref{lem:T}, and also (\ref{eq:xi0isximistaum}) remain valid in this
case.
We thus have the following.

\begin{theorem} \label{thm:LGm1ACr}
For any oriented link $L$, there holds:
\begin{eqnarray*}
LG^{m,1}_L(\tau,e^{\pi\sqrt{-1}\,r/m}) =\Delta_L(\tau^{2m}),
\end{eqnarray*}
where $r$ is any integer such that $r$ and $m$ are relatively prime.
\end{theorem}

In particular, via the choice $r=-1$ and the use of
symmetry~(\ref{eq:SymmetryofLGmn}), we immediately deduce the
following.
\begin{corollary} \label{cor:LG1n}
For any oriented link $L$, there holds:
\begin{eqnarray*}
LG^{1,m}_L(\tau,e^{\pi\sqrt{-1}/m}) =\Delta_L(\tau^{2m}).
\end{eqnarray*}
\end{corollary}

\section{Extensions}

To conclude, we believe that Theorem~\ref{thm:LGm1AC} can be extended to
a similar statement for $LG^{m,n}$.

\begin{conjecture} \label{conj:LGmnAC}
For any oriented link $L$, there holds:
\begin{eqnarray}
LG^{m,n}_L(\tau,e^{\pi\sqrt{-1}/m}) =\Delta_L(\tau^{2m})^n,
\label{eq:LGmnAC}
\end{eqnarray}
and equivalently, by the symmetry~(\ref{eq:SymmetryofLGmn}):
\begin{eqnarray}
LG^{m,n}_L(\tau,e^{\pi\sqrt{-1}/n}) =\Delta_L(\tau^{2n})^m.
\label{eq:LGmnACprime}
\end{eqnarray}
\end{conjecture}
Thus, for a given invariant $LG^{m,n}$, there are two distinct
reductions which recover $\Delta$; note that Theorem~\ref{thm:LGm1AC}
and Corollary~\ref{cor:LG1n} are particular cases of
Conjecture~\ref{conj:LGmnAC}.
We mention that the considerations leading to Theorem~\ref{thm:LGm1ACr}
also lead to the obvious generalisation of
Conjecture~\ref{conj:LGmnAC}.

These relations are initially surprising in that neither is symmetric in
$m$ and $n$; however, we have a range of evidence to support them.
For instance, we can verify~(\ref{eq:LGmnACprime}) for $LG^{2,1}$ for
closed $2$-braids $\widehat{\sigma^k}$, where $\sigma$ is the generator
for the braid group $B_2$.
To that end, with reference to~(\ref{eq:R}), we have for $LG^{2,1}$:
\begin{eqnarray*}
R =q^{-2\alpha}P_0 -P_1 +q^{2\alpha+2}P_2
  =\tau^2P_0 -P_1 +\tau^{-2}q^2P_2,
\end{eqnarray*}
thus:
\begin{eqnarray*}
R^k =\tau^{2k}P_0 +(-1)^kP_1 +\tau^{-2k}q^{2k}P_2.
\end{eqnarray*}
Specialising the formulae of Lemma 2, we have:
\begin{eqnarray*}
\mathrm{cl}(P_0) &=&
\frac{\tau -\tau^{-1}}
     {(\tau q +\tau^{-1} q^{-1})(\tau^2 q -\tau^{-2} q^{-1})}
\, \mathrm{id}_V, \\
\mathrm{cl}(P_1) &=&
\frac{-(q+q^{-1})}
     {(\tau +\tau^{-1})(\tau q +\tau^{-1} q^{-1})}
\, \mathrm{id}_V, \\
\mathrm{cl}(P_2) &=&
\frac{\tau q -\tau^{-1} q^{-1}}
     {(\tau + \tau^{-1})(\tau^2 q - \tau^{-2} q^{-1})}
\, \mathrm{id}_V.
\end{eqnarray*}
So, in the substitution $q=-1$, denoting $X|_{q=-1}$ by
$\overline{\overline{X}}$:
\begin{eqnarray*}
           \overline{\overline{\mathrm{cl}(R^k)}}
= \tau^{2k}\overline{\overline{\mathrm{cl}(P_0)}}
    +(-1)^k\overline{\overline{\mathrm{cl}(P_1)}}
+\tau^{-2k}\overline{\overline{\mathrm{cl}(P_2)}},
\end{eqnarray*}
where
$\overline{\overline{\mathrm{cl}(P_0)}}
=-{\textstyle \frac{1}{2}}
 \overline{\overline{\mathrm{cl}(P_1)}}
=\overline{\overline{\mathrm{cl}(P_2)}}
=(\tau + \tau^{-1})^{-2} \, \mathrm{id}_{\overline{\overline{V}}}$,
and so:
\begin{eqnarray*}
LG^{2,1}_{\widehat{\sigma^k}}(\tau,-1)
=\left( \frac{\tau^{k} -(-\tau)^{-k}}{\tau +\tau^{-1}} \right)^2.
\end{eqnarray*}
Then, for $\Delta(\tau^2)\equiv LG^{1,1}(\tau)$, we have
$R =\tau P_0 -\tau^{-1} P_1$, where
\begin{eqnarray*}
\mathrm{cl}(P_0) =-\mathrm{cl}(P_1)
=(\tau + \tau^{-1})^{-1} \, \mathrm{id}_V.
\end{eqnarray*}
Hence
$\Delta_{\widehat{\sigma^k}}(\tau^2)
=({\tau^{k}-(-\tau)^{-k}})/({\tau+\tau^{-1}})$,
and thus:
\begin{eqnarray*}
LG^{2,1}_{\widehat{\sigma^k}}(\tau,-1)
=\Delta_{\widehat{\sigma^k}}(\tau^2)^2.
\end{eqnarray*}
Similarly, we can verify Conjecture~\ref{conj:LGmnAC} for $LG^{2,2}$ for
closed $2$-braids $\widehat{\sigma^k}$, using formulae derived
in~\cite{GouldLinksZhang96} (specifically, formula (71) and explicit
details described in later sections).
That is, we have:
\begin{eqnarray*}
LG^{2,2}_{\widehat{\sigma^k}}(\tau,e^{\pi\sqrt{-1}/2})
= \Delta_{\widehat{\sigma^k}}(\tau^{4})^2.
\end{eqnarray*}
Lastly, we have also been able to computationally verify
(\ref{eq:LGmnACprime}) for $LG^{2,1}$ for a range of prime knots using
the state model method of evaluation for $LG^{2,1}$ described
in~\cite{DeWit00}.
Specifically, this has been done for a selection of $4310$ prime knots
of up to $14$ crossings, including all prime knots of up to $10$
crossings.
Beyond that, using formula (71) of~\cite{GouldLinksZhang96}, we have
also verified that~(\ref{eq:LGmnAC}) holds for $LG^{m,n}$ for all
$m,n\leqslant 5$, for closed $2$-braids $\widehat{\sigma^k}$ for
$k=2,\dots,6$ (and thereby, for all $0\leqslant|k|\leqslant 6$).

\vspace{\baselineskip}

Now, let $GLZ^n$ denote the invariants proposed
in~\cite{GouldLinksZhang96} associated with the $U_q[osp(2|2n)]$
superalgebras.
We can state a similar result to Conjecture~\ref{conj:LGmnAC}.

\begin{conjecture}
For any oriented link $L$, there holds
\begin{eqnarray*}
GLZ^n_L(\tau,e^{\pi\sqrt{-1}/2}) =\Delta_L(\tau^{4})^n.
\end{eqnarray*}
\end{conjecture}

As $osp(2|2)\cong sl(2|1)$, we have $GLZ^1\equiv LG^{2,1}$, so this
conjecture is true for $n=1$.
Further evidence for it is that via similar considerations to the above
using results from~\cite{GouldLinksZhang96}, we have confirmed that:
\begin{eqnarray*}
GLZ^2_{\widehat{\sigma^k}}(\tau,e^{\pi\sqrt{-1}/2})
=\Delta_{\widehat{\sigma^k}}(\tau^{4})^2.
\end{eqnarray*}
The difficulty in proving these conjectures lies in the fact that in
general $\Delta(t)^n$ satisfies higher-order skein relations.
One could begin by establishing that $LG^{m,n}$ and $GLZ^n$, at the
appropriate values of $q$, satisfy the same skein relations as
$\Delta(t)^n$.
For example, for $LG^{2,1}$, two such skein relations are known, and
these may be used to evaluate the invariant for all algebraic
links~\cite{Ishii03skein}.
For $q=-1$, we have checked that both skein relations reduce to ones
which are satisfied by $\Delta(t)^2$, which confirms
that~(\ref{eq:LGmnACprime}) holds for $LG^{2,1}$ for a vast class of
links.
However, it is not clear that these two skein relations are sufficient
to determine $LG^{2,1}$ for any arbitrary link.
More generally, for either $LG^{m,n}$ or $GLZ^n$, the only
easily-determined skein relation is that corresponding to the
characteristic identity satisfied by $R$ (illustrated for the $LG^{m,1}$
case in~(\ref{eq:characteristicidentityofRm1})).
Additional skein relations are generally not known.

\subsection*{Acknowledgements}

We thank Michael Gagen and Mark Gould for helpful comments.
Jon Links thanks the Australian Research Council for support through an
Australian Research Fellowship.

\Addresses\recd

\end{document}